\documentstyle[epsfig]{article}

\newtheorem{prop}[equation]{Proposition}
\newtheorem{cor}[equation]{Corollary}
\newtheorem{theorem}[equation]{Theorem}
\newtheorem{lemma}[equation]{Lemma}
\newtheorem{remar}[equation]{Remark}
\newtheorem{definitio}[equation]{Definition}

\newtheorem{notat}[equation]{Notations}
\newtheorem{nott}[equation]{Notation}
\newenvironment{remark}{\begin{remar} \rm }{\end{remar}}
\newenvironment{defin}{\begin{definitio} \rm }{\end{definitio}}
\newenvironment{nota}{\begin{notat} \rm }{\end{notat}}
\newenvironment{note}{\begin{nott} \rm }{\end{nott}}

\catcode`\@=11
\def\st{\mathop {\operator@font st}\nolimits}
\def\eqalign#1{\null \,\vcenter {\openup \jot \m@th \ialign 
{\strut \hfil $\displaystyle {##}$&$\displaystyle {{}##}$\hfil \crcr 
#1\crcr }}\,}\catcode`\@=12

\font\tensym=msbm10    
\font\sevensym=msbm7
\font\fivesym=msbm5
\newfam\symfam
\textfont\symfam=\tensym
\scriptfont\symfam=\sevensym
\scriptscriptfont\symfam=\fivesym
\def\sym{\fam\symfam\tensym}

\def\R{{\sym R}}
\def\N{{\sym N}}
\def\Z{{\sym Z}}

\def\qs{\forall\,}
\def\imp{\ \hbox{$\Relbar\mkern -10mu\Rightarrow$~}}
\def\b{$\bullet$\ }
\let\0=\emptyset
\def\diagram#1{\def\normalbaselines {\lineskip=7pt\baselineskip=0pt
\lineskiplimit=1pt}\matrix{#1}}
\def\hfl#1{\smash{\mathop{\hbox to 10mm{\rightarrowfill}}\limits^{\textstyle
#1}}}

\def\vf{\left\downarrow\vbox to 5mm{}\right.}

\let\Phi=\Psi
\input epsf
\def\A{{\cal A}}
\def\g{\Gamma}
\def\od{G^3}
\def\odn{\od_n}
\def\B{F}
\def\dl{D_{n,k}}
\def\ed{{\cal D}}
\def\edn{\ed_n}
\def\ew{{\cal D}^W}
\def\ewn{{\cal D}^W_n}
\def\dw{D^W_{n,k}}
\def\eg{{\cal D}^g}
\let\gu=\gamma
\def\egn{{\eg_n}}
\def\dg{D^g_{n,k}}
\def\cem{C_\g^K}
\def\c{{C_\g}}
\def\cc{{C'_\g}}
\def\wg{{\cal W}_\ph(\gu)}
\def\ce{{\cal E}}
\def\wp{{\cal W}_{\ph'}(\gu)}
\def\wx{W_x}
\let \inc=\subseteq
\let\O=\Omega
\def\h{{\cal H}}
\let\ph=\varphi
\def\cb{\overline \c}
\def\gr{\overline C(G)}
\def\f{f}
\def\a{\alpha}\def\ss{\co_S}
\def\Si{\Sigma}
\def\sd{S^2}
\let \ps=\psi\def\rt{\R^3}
\def\dd{D^2}
\def\an{a_n}
\def\ano{a}
\def\anp{a_{n'}}
\def\D{l_x}
\def\lo{l_0}
\def\al{Q}

\def\s{\sigma} 
\let\om=\omega
\let\th=\theta
\def\bb{{\bar B}}
\def\enk{\ed_{n,k-1}}
\def\ba{{\bar A}}
\def\vep#1{\vcenter{\epsfbox{#1.eps}}}
\def\epfb#1{\left(\,\vcenter{\epsfbox{#1.eps}}\,\right)}
\let\de=\Delta
\let\tens=\otimes
\def\bu{{\bar U}}
\def\ha{{\widehat A}}
\def\ra{R_A}
\def\hh{\h'}
\def\rb{R_\ha}
\def\1{+\infty}
\def\cqfd{\unskip\kern 6pt\penalty 500\raise -2pt
 \hbox{\vrule\vbox to8pt{\hrule width 4pt\vfill\hrule}\vrule}\par}

\let\lra=\longrightarrow
\def\co{C^0}

\title{Rationality Results for the Configuration Space Integral of Knots}
\author{Sylvain Poirier}
\begin{document}
\maketitle

{\def\thefootnote{\relax}

\begin{abstract}
The perturbative Chern-Simons theory for knots in 
Euclidean space is a linear combination of integrals on configuration spaces.
This has been successively studied by Bott and Taubes, Altschuler and
Freidel, and Yang. We study it again in terms of degree theory,
with a new choice of compactification. This paper is self-contained 
%(excepts for the algebra of diagrams), 
and proves some old and new results, especially a rationality result
with some information on the denominators.
\end{abstract}
\tableofcontents
\section{Introduction}
\setcounter{equation}{0}
\subsection{The Hopf algebra $\A$ of diagrams}

A {\it diagram} can be defined in two ways:

First, it can be seen as the two following data.

\b A connected trivalent graph, where a trivalent graph is a finite set 
$V$ of vertices together with a finite
set $E$ of edges and a relation between them, such that every edge
is related to two vertices and every vertex belongs to precisely three edges.

\b A preferred oriented cycle (made of vertices and edges alternatively), 
and an orientation on every vertex which
does not belong to this cycle. (An orientation is a cyclic order on the set of the
three edges containing this vertex). The cycle
must be nonempty. Thus it contains at least two vertices if the
graph is nonempty.

For the second point of view, which is the one we use in this article,
the edges which belong to the cycle are deleted, and replaced by a cyclic
order on its vertices which are now univalent. So this graph can be 
disconnected, but each connected component has univalent vertices.

\begin{nota}
Let $\od$ denote the set of these diagrams up to isomorphism.

The {\it degree} of a diagram is half its number of vertices.
The set of degree $n$ diagrams is denoted by $\odn$.

We denote by $\A_n$ the real vector space generated by $\odn$ and quotiented
by the following AS, IHX and STU relations:

The relation AS says that if $\g$ and $\g'$ differ by the 
orientation of one vertex, then $\g'=-\g$.

The IHX and STU relations are respectively defined by the following
formulas. 
$$\vep{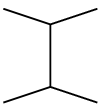}=\vep{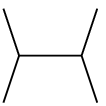}-\vep{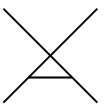}\qquad,\qquad\vep{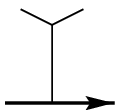}=\vep{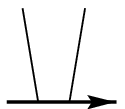}-\vep{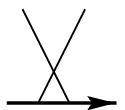}$$
These formulas relate diagrams which are identical outside
one place, where they differ according to the figures.
By convention, the orientation at each vertex is the given counterclockwise
order of the edges.

The image of a diagram $\g\in\odn$ in
the space $\A_n$ is denoted by $[\g]$.

The space $\A_0$ is one-dimensional. It is generated by the empty diagram.

The space $\A_1$ is one-dimensional. It is generated by a single diagram
denoted by $\th$, for it looks like the letter $\th$.

Now let $\A$ denote the space $$\A=\prod_{n\in\N}\A_n.$$
This space has the structure of a graded Hopf algebra (as in \cite{mm}).
\end{nota}

The algebra structure we take here is the connected sum of diagrams,
obtained by cutting each cycle at some place and gluing them together
(and the result is independent on the places where the cyclic orders 
are cut, modulo the AS, IHX and STU relations: see \cite{bn}).

We recall the coalgebra structure: for any diagram $\g$, let $X$ be the set
of its connected components. Then
$$\de(\g)=\sum_{A\inc X} \g_A \tens\g_{X-A}$$
where $\g_A$ is the diagram obtained by deleting from $\g$ all components
that do not belong to $A$.

An element $x$ of a Hopf algebra is said to be {\it primitive} iff
$$\de(x)=x\tens1+1\tens x.$$
The set of primitive elements of $\A_n$ is the space generated by the connected
degree $n$ diagrams.

An element $x$ of a Hopf algebra is said to be {\it grouplike} iff $x\not=0$
and
$$\de(x)=x\tens x.$$
It is well-known that the grouplike elements of $\A$ are the exponentials
of its primitive elements.

In the following we denote by $\ed$ a system of arbitrary representants of
the diagram classes in $\od$ which only differ by the orientations
of trivalent vertices.

\subsection{Statement of the main theorem}

Let $K$ be a knot, that is, a smooth imbedding: $S^1\hookrightarrow \R^3$. 
%The constructions discussed below could be naturally generalised to links; we 
%shall generalise them to braids and tangles in the second part of this article.

Our purpose will be to study the following expression depending on the knot $K$
with values in $\A$. (It represents the perturbative expression of
the Chern-Simons theory applied to $K$):
\begin{equation}
\label{feyn}
Z(K)=1+\mathop{\sum_{\g\in\ed}}\limits_{\deg \g\geq1}{I_K(\g)\over|\g|} [\g],
\end{equation}
where: 

$[\g]$ is the image of $\g$ in $\A$;

$|\g|$ is the number of automorphisms of $\g$ (considered without 
vertex-orientations); and
$$I_K(\g)=\pm\int_\c\bigwedge_{e\;\rm edge}\phi_e^*\om$$
where $\c$ is the configuration space of imbeddings of the set of vertices
of $\g$ in $\R^3$ in which the restriction to univalent vertices factorises
through $K$, respecting the cyclic order; $\phi_e$ is a map from $\c$ to
$S^2$ defined by the vector connecting the two ends of the edge $e$ in
$\R^3$, and $ \om$ is the standard volume form on $S^2$ with total mass 1.
The sign of that expression is not well-defined because the space $\c$ is
not oriented yet. This will be precised later.

This expression has been studied by Bott and Taubes \cite{bt}, who
first proved the convergence of the integrals $I_K(\g)$ by compactifying the
configuration space $\c$ into a manifold with corners $\cb$ on which the
above differential forms extend smoothly.

Then, Bott and Taubes applied the Stokes theorem to express the variations
of $I_K(\g)$ during an isotopy of the knot, in terms of integrals on
the faces of $\cb$. They found equalities between the integrals
at certain faces of different $\cb$, where the diagrams $\g$ are the terms of
an IHX or STU relation, and proved cancellation on the other faces except
those of a special kind. They called ``anomaly'' the contribution
of these last faces. Its global effect can be expressed by a primitive element
$a\in \A$ which is independent of the knot.

The degree one part of $a$ is $[\th]$. Yang \cite{ya} claims that the terms of
other degrees cancel (in other words, $a=[\th]$). Anyway,
Altschuler and Freidel \cite{af} have shown (and we shall prove 
it in another way with Proposition \ref{group}), that
$Z$ is grouplike and that the variations of $Z$ on an isotopy class of 
knots can be expressed
by the formula $$Z=Z_0 \exp\left({I_K(\th)\over2}\ano\right)$$
where $Z_0$ only depends on the isotopy class of the knot. The integral
$I_K(\th)$ is the well-known Gauss integral, which takes all real values
on each isotopy class of knots. They also prove that $Z_0$ is a universal 
Vassiliev invariant (and we shall not do it again).

The aim of this article will be to prove the following theorem, after
the notations.

\begin{nota}

Let $\th$ be the diagram with only one edge.

For $n\geq1$, let $Z_n(K)$ be the degree $n$ part of $Z(K)$.

For a diagram $\g$, let $u_\g$ be the number of univalent vertices of $\g$.

Let $k$ be an integer, and let $\A_n^k$ be the quotient space of 
$\A_n$ by the subspace generated by the diagrams $\g$ such that $u_\g<k$,
and each connected component of $\g$ is either a chord 
(with no trivalent vertex) or has at least three univalent vertices.

Let $Z_n^k$ be the image of $Z_n$ in $\A_n^k$.
\end{nota}

Note that $\A_1^k=\A_1$ if $k\leq 2$, and that $\A_n^k=\A_n$ if $n>1$ and
$k\leq3$.

\begin{theorem}
\label{mult}
We suppose that the knot $K$ is such that $I_K(\th)$ is an integer. 
Let $n\geq1$ and $k\in\Z$. Then $Z_n^k$ belongs to the integral lattice
generated by the elements$${(u_\g-k)!\over (3n-k)!\;2^{3n-u_\g}}[\g]$$
in the space $\A_n^k$, where $\g$ runs over the
diagrams such that $u_\g\geq k$ and each non-simply connected 
component of $\g$ has at least four 
univalent vertices.
%and if for each $\g$, $F(\g)$ is a multiple of
%$${(3n-k)!\over (u_\g-k)!}2^{3n-u_\g},$$
%then $F(Z_n(K))$ is an integer.
\end{theorem}

This is the first in a series of two articles which aims to prove
(with the help of the zero-anomaly result $a=[\th]$ of Yang) the
equality between the integral $Z$ above and the Kontsevich integral. 
There is also a rationality result for 
the Kontsevich integral by \cite{le}, which we recall:

The degree $n$ part of the Kontsevich integral is a linear combination
of chord diagrams (or of all diagrams, which is the same), where the
coefficients are multiples of$$1\over (1!2!\cdots n!)^4(n+1)$$

For large enough values of $n$, these denominators are greater than
those of the above theorem, but have smaller prime factors.

\bigskip
I would like to express my thanks to Christine Lescop for the help
during the preparation of this paper. 
I also thank Gregor Masbaum, Pierre Vogel, Su-win Yang, Dylan Thurston, 
Dror Bar-Natan, and Michael Polyak for their advice and comments.

\section{Sketch of proof of Theorem \ref{mult}}
\label{def}
\setcounter{equation}{0}

The constructions we shall make in this article will be based
on labellings and orientations of the edges. We shall call that ``labellings"
of diagrams.
We are going to write this on the very definitions of the diagrams.

We also fix the degree of diagrams to be equal to $n$, but since the number of
edges is not constant and since we must anyway fix the set of their labels,
there will be some unused labels which will be called ``absent edges".
Moreover, we shall only work with the diagrams which verify a special
condition (being ``triply connected" to $U$), and which have a number of 
edges lower than or equal to the number $N$ of labels, and thus a number
of univalent vertices $u_\g\geq k=3n-N$.
% of class $C^1$, and $C^2$ by parts.
%For a graph $\g$, let $u_\g$ denote the number of univalent vertices of $\g$.
%Let $n\geq 2,\ k\geq 3$ be fixed integers, and $F$ be a weight function in degree $n$,

\begin{nota}
Let $K$ be a fixed knot; let $n\geq 1,\ k\leq2n$ be fixed integers and 
$N=3n-k$.

Let $E=\{e_1,\cdots, e_N\}$ be the set of labelled ``edges" 
$e_i=\{2i-1,2i\}$.
The set $X=\{1,\cdots, 2N\}$ is called the set of half-edges.

If $A$ is a set of sets, $\cup A$ will denote the union of 
elements of $A$.
\end{nota}

\begin{defin}\label{diag}%\proclaim Definition. 
We shall denote $\dl$ the set of 4-uples 
$\g=(U,\s,T,E^v)$,
called {\it labelled diagrams}, 
%(i.e. graphs with labelled and oriented edges), 
which respect the following conditions:

\b $U$ is the set of univalent vertices: it is a set of singletons included
in $X$, and $\s$ is a circular
permutation of $U$. We define $u_\g=\#U$.

\b $E^v\inc E$ is the set of ``visible edges", and labels the edges of 
the diagram. We denote $E^a=E-E^v$, the set of ``absent edges", 

\b $T$ is the set of trivalent vertices, which are parts of $X$ with 
three elements belonging to three different edges.

\b The set $V=U\cup T$ of all vertices of $\g$ is a partition of $\cup E^v$; in
other words, $E^a\cup T\cup  U$ is a partition of $X$.

\b We impose 
$\#V=2n$ (which is equivalent to $u_\g=\#E^a+k$ by the above conditions).

\b $\g$ must be triply connected to $U$, with respect to the following 
definition 2.5.
\end{defin}

\begin{nota}%\proclaim Notation. 
Let $A$ be any nonempty part of $V=U\cup T$. We denote
$$\eqalign{E_A&=\{e\in E|e\inc\cup A\}\cr E'_A&=\{e\in E|\#(e\cap\cup A)=1\}}.$$
The cardinalities of these sets are related by the formula
\begin{equation}
\label{count}
2\#E_A+\#E'_A=3\#(A\cap T)+\#(A\cap U)\quad(=\#\cup A)
\end{equation}
\end{nota}

\begin{defin} 
A diagram $\g$ is said to be {\it triply connected to $U$} if for any 
$A\inc T,\ \#A>1$ we have $\#E'_A\geq 3$.
\end{defin}

%Note this implies that the two ends of an edge belong to different vertices.
Note that the assumption $k\leq 2n$ is equivalent to the fact that
$\dl\not=\0$, and implies that $\dl$ contains all chord diagrams (diagrams
with $T=\0$).

\begin{defin} An {\it isomorphism} between two diagrams 
$\g=(U_\g,\s_\g,T_\g,E^a_\g)$ and $\g'=\-(U_{\g'},\s_{\g'},T_{\g'},E_{\g'}^a)$ 
(or {\it change of labelling}) is
a bijection between $\cup E^v_\g$ and $\cup E_{\g'}^v$ which carries
$E^v_\g$ to $E^v_{\g'}$, $T_\g$ to $T_{\g'}$, $U_\g$ to $U_{\g'}$ and 
$\s_\g$ to $\s_{\g'}$.
\end{defin}

We shall now introduce the notion of configuration space.
\subsection{The configuration spaces $\c$ and $\cc$ and the maps 
$\Phi_\g$ and $\Phi'_\g$}

\begin{defin} The configuration space $C_\g$ of a diagram $\g$
(for the knot $K$) is the space of maps $\f$ from $V(\g)$ to $\R^3$
such that :

1) Any two vertices connected by an edge are mapped to different points
(this is slightly better than to suppose all the images different and gives the
same integral). %This condition will be implicit for all the other 
%``configuration spaces" defined in this article.

2) $\f_{|U}$ is the composition of $K$ with an injection of $U$ into $S^1$ 
which respects the cyclic order defined by $\s$.
\end{defin}%restriction ???

The dimension of $\c$ is $\#U+3\#T=2\#E^v$.
Denoting also by $\f$ the induced map from $\cup E^v$ to $\R^3$, 
we have a canonical map: $$\eqalign{\Phi_\g:C_\g&\longrightarrow (S^2)^{E^v}\cr 
f&\longmapsto 
\left({\f(2i)-\f(2i-1)\over|\f(2i)-\f(2i-1)|}\right)_{e_i\in E^v}.
\cr}$$

This map extends to a map $\Phi'_\g$ from the space $\cc=\c \times (S^2)^{E^a}$
to $(S^2)^E$, which is the product of $\Phi$ by the identity on 
$(S^2)^{E^a}$: if $\f\in\cc$, $\Phi'(\f)(e)=\Phi(f_{|V})(e)$ 
if $e\in E^v$ and $\Phi'(\f)(e)=f(e)$ if $e\in E^a$.

\subsection{Orientation of $\cc$}\label{ori}
An {\bf orientation} $o$ of a diagram $\g$ is the datum of an orientation of
all its trivalent vertices. A diagram $\g$ is said to be {\bf oriented}
if it is given an orientation. We shall now define an orientation of the 
space $\cc$ depending on the orientation $o$ of $\g$. It will be done by
labelling all components of a natural local coordinate system of $\cc$ 
by the elements of $X$. These components are:

\b One component in the oriented circle $S^1$ for each univalent vertex.

\b the three standard components in $\rt$ for each trivalent vertex.

\b For each absent edge we have the two components of an oriented local 
coordinate
system in $S^2$, with respect to the standard orientation of $\sd$.

Since a univalent vertex is a singleton of $X$, its coordinate on $S^1$
is already labelled. The two components of each $e_i\in E^a$
are labelled by $2i-1$ and $2i$, respectively.% label components of a local
%coordinate system for the point $\f(a_i)\in S^2$ with its standard 
%orientation.

%First, the elements of $U$ canonically label the components of $\f_{|U}$, 
%each component respecting the orientation of $S^1$.

As for the trivalent vertices, let us represent
$o$ by a family of bijections $(h_t)_{t\in T}$ from each trivalent vertex $t$
to the set $\{1,2,3\}$.
Then, for all $x\in t$, give the label $x$ to the $h(x)$-th component of 
$\f(t)$ in $\R^3$.

Thus we have an orientation of $\cc$, which only depends on $o$
and is reversed when we change the orientation of one vertex.

For a labelled oriented diagram $(\g,o)$, we can now define the integral 
$$I_K(\g,o)=\int_\cc {\Phi'_\g}^{\kern -2pt*}\O=\int_\c \Phi_\g^*\O_{E^v}$$
where $\O$ is the standard volume form on $(S^2)^E$ with total mass 1,
and $\O_{E^v}$ is the one of $(S^2)^{E^v}$.

The convergence of this integral results from the existence of a 
compactification of $\c$ and thus of $\cc$ to a manifold with corners,
on which the map $\Phi_\g$ extends smoothly. 

%This could be generalised to the case where the knot $K$ is of class $C^2$,
%with a map $\Phi$ of class $C^1$ near the boundary (and of class $C^2$ outside
%it). It could even be generalised to a knot of class $C^1$ and $C^2$ by 
%closed parts, but then the boundary would be more complicated.

Note that we have $I(\g,-o)=-I(\g,o)$ and that $I(\g,o)$ is independent 
of the labelling of $\g$: exchanging two edges preserves both 
orientations of $\cc$ and $(S^2)^E$, and reversing one edge changes both 
orientations (in other words, it changes both the orientation of
$\cc$ and the sign of the form ${\Phi'_\g}^*\O$), because the antipodal 
map reverses the orientation of $S^2$.

So, the product $I_K(\g)[\g]$ is well-defined. It does not depend 
on the orientation $o$ of $\g$.

In all the following, we fix an arbitrary choice of orientations 
of the diagrams in $\dl$, and it will be used implicitly,
with no consequence on the final results.

\subsection{The map $\Phi$ and its degree}\label{phi}

Let $\a $ be a map from $\dl$ to a $\Z$-module $M$. Denote
by $\Phi$ the union of the maps $\Phi'_\g$ on the different $\cc$ for
$\g\in \dl$.

For a generic $Q\in(\sd)^E$, consider the expression 
$$\sum_{\g\in\dl}\a (\g)n_\g(Q)$$ 
with values in $M$, where $n_\g(Q)$ is the number of preimages $Q'$  of 
$Q$ by the map $\Phi'_\g$, with $Q'$ counted positively if the
differential d$\Phi$ at $Q'$ preserves the orientation, and negatively
if it reverses the orientation.

According to the degree theory, this expression is locally constant with
respect to $Q$, except on the images by $\Phi$ of the boundaries of
the spaces $\cb$ that will be obtained by compactifying $\cc$.
%for each $\g$, the corresponding term produces a difference of $\a (\g)$
%between the sides of $\Phi(\partial \cb)$.

But if the expression was constant with
respect to $Q$, it would define the degree of the map $\Phi$
(which depends on $\a $)
on the formal linear combination of the spaces $\cc$:
$$C=\sum_{\g\in \dl}\a (\g) \cc$$
$$\deg\Phi=\sum_{\g\in \dl}\a (\g)n_\g(Q)$$%\quad\qs Q\in(\sd)^E.$$
for any generic $Q$. This would work that way if we could find that 
some of the faces of the $\cb$
have images by $\Phi $ with codimension at least 2,
and that the other faces can be grouped into equivalence classes which verify
the two conditions: first, two equivalent faces of spaces $\cb$ and $\overline
C_{\g'}$ can be identified with each other by an identification which 
carries $\Phi'_\g$ to $\Phi'_{\g'}$, restricted to these faces.
%First, they identify all together such that the map $\Phi$ is
%the same on them. 
Second, the sum (with signs...) of the coefficients $\a $ of the diagrams 
of the faces in a class cancel in $M$ (this is the ``gluing condition'').

In particular, to prove Theorem \ref{mult},
we will consider the map $\a $ from $\dl$ to $\A_n^k$ defined by
\begin{equation}\label{coef}
\a (\g)={(u_\g-k)!\over (3n-k)!\;2^{3n-u_\g}}[\g].
\end{equation}
where $[\g]$ denotes the class of $\g$ in the space $\A_n^k$.

If the above cancellation of the face contributions worked, 
we would deduce that the element of $\A_n^k$ 
$$\sum_{\g\in \dl}\a (\g) I_K(\g)=\sum_{\g\in \dl}\a (\g) 
\int_\cc {\Phi'_\g}^{\kern -2pt*}\O=\deg\Phi$$
is in the lattice generated by the $\a (\g)$ such that 
$\Phi(\cc)$ has interior points.

Unfortunately we cannot make the faces of the $\cb$ cancel each other 
in general. But we will find extra spaces to add to the combination 
in order to obtain such a cancellation.

Anyway, let us already explain roughly how Theorem \ref{mult} is deduced 
from this construction.
We are going to prove the three following facts:

1) The integral on the new spaces can be chosen to vanish if $I_K(\th)$
is an integer.

2) The above gluing conditions are satisfied for the expression (\ref{coef}).

3) For any $\g\in\edn$, if $\g$ is not represented in $\dl$, then
$I_K(\g) [\g]=0$.

Then, from Equation (\ref{coef}) (where $\a $ is constant on the
isomorphism classes), we deduce
\begin{equation}
\label{som}
\qs \g_0\in\dl,\quad {[\g_0]\over|\g_0|}=\sum_{\g\in \dl,\;\g\sim\g_0}\a (\g)
\end{equation}

Indeed, we have $3n-u_\g=\#E^v$, and the number of diagrams 
which are isomorphic to a given diagram $\g$ is the product of the number
$(3n-k)!\over(u_\g-k)!$ of injections from the set of edges 
of a diagram $\g$ to the set $E$, by the number $2^{3n-u_\g}$
of orientations of the edges, divided by the number $|\g|$ 
of automorphisms of $\g$.

Then, by Fact 3),
$$\sum_{\g\in \dl}\a (\g) I_K(\g)=\sum_{\g\in\edn}{I_K(\g)\over|\g|} 
[\g],$$in other words, $\deg \Phi=Z_n^k$. 

The generators of the integral lattice
containing deg $\Phi$ are the $\a (\g)$ such that
$[\g]\not=0$ in $\A_n^k$ and $\Phi(\cc)$ has interior points. There remains
to show these generators satisfy the conditions stated in Theorem \ref{mult}.
This is the consequence of the following lemma, 
which will be proved in the next section:

\begin{lemma}
\label{quad}
If $\Phi_\g(\c)$ has interior points, each connected
component of $\g$ which is non-simply connected has at least four 
univalent vertices.
\end{lemma}

So if $[\g]\not=0$ in $\A_n^k$ and $\Phi(\c)$ has interior points,
then $u_\g\geq k$, %il s'ensuit que $[\g]=0$ dans $\A_n^k$, 
by definition of $\A_n^k$.

Let us now prove Fact 3). The proofs of Facts 1) and 2) will
be obtained in Section 6.
\medskip

The diagrams $\g\in \edn$ which are represented in $\dl$ 
are those triply connected to $U$ and such that
$u_\g\geq k$. So if a diagram $\g\in\edn$ is not represented in $\dl$,
there are two cases:

First case: $\g$ is not triply connected to $U$. Then we deduce
$I_K(\g)=0$ from the following lemma, which will be fully used in Section $5$
(here a codimension 1 is enough).

\begin{lemma}
\label{trip}
If $\g$ is not triply connected to $U$ then $\Phi(\cc)$ lies in a 
space of codimension 2 in $(S^2)^E  $.
\end{lemma}

Second case: $\g$ is triply connected to $U$ and $u_\g<k$: this
implies $[\g]=0$ (in $\A_n^k$).

\section{Compactification}
\setcounter{equation}{0}

\begin{defin} If $A$ is a finite set with at least two elements, 
$C^A$ denotes the space $(\R^3)^A/TD$ of nonconstant maps from 
$A$ to $\R^3$ quotiented by the translation-dilations group 
(that is the group of translations and positive homotheties of $\R^3$).
\end{defin}

Note that the space $C^A$ is diffeomorphic to the sphere $S^{3\#A-4}$.
Thus it is compact. 

\medskip
\noindent {\sc Proof of lemmas \ref{quad} and \ref{trip}: }Let 
$A\inc T,\ \#A>1$. 

Consider the following commutative diagram:
$$\diagram{\c&\hfl{\Phi_\g}&(S^2)^{E^v}\cr\vf&&\vf\cr C^A&\hfl{}&(S^2)^{E_A}}$$

According to (\ref{count}) we have 
$$2\#E_A+\#E'_A=3\#A%\leq2\#E_A+2
\ \imp\ \dim C^A=\dim(S^2)^{E_A}+\#E'_A-4.$$

For Lemma \ref{quad}, if $\g$ has a connected component $B$ which is
not simply connected and which has less than 4 univalent vertices, 
then by taking $A=T\cap B$ 
we have $\#A>1$, and $\#E'_A\leq3$.

For Lemma \ref{trip}, if $\g$ is not triply connected to $U$, there
exists $A\inc T,\ \#A>1$ such that $\#E'_A\leq 2$, so the image
of $\Phi_\g$ has codimension at least 2.\cqfd

\medskip

Now, before defining a compactification of our space $\c$, we shall 
describe a general compactification tool (which will be of much use in
our next article). 

\subsection{Compactification of a configuration space $\co(G)$ of a graph}

\begin{nota} Let $G$ be a graph defined as a couple $G=(V,E)$ 
where $V$ is the set of ``vertices'' and $E$
is the set of ``edges'': $V$ is a finite set and $E$ is a set of 
pairs of elements of $V$. We suppose
that $\#V\geq 2$ and that $G$ is connected.

We define the configuration space $\co(G)$ as the space of maps
from $V$ to $\R^3$ which map the two ends of every edge to two different points, modulo
the dilations and translations (it is a dense open part of $C^V$).

By an abuse of notation, any part $A$ of $V$ will be considered as the graph
with vertices $A$ and edges $E_A$ (the set of edges $e\in E$ such that $e\inc
\cup A$).  

Let $R$ be the set of parts $A$ of $V$ such that 
$\#A\geq 2$ and $A$ is connected;

$\h=\prod_{A\in R}C^A$.
\end{nota}

Note that the manifold $\h$ is compact.

We are going to define a 
canonical imbedding of $\co(G)$ in $\h$, and define the compactification
$\gr$ of $\co(G)$ to be its closure in $\h$. 
Then, we are going to see that $\gr$ is a manifold
with corners, and describe completely all its strata as abstract spaces
and their imbeddings in $\h$. 

But first, let us see the simplest case of the imbedding of $\co(G)$:

This map is defined by restricting any $f\in \co(G)$ 
to each of the $A\in R$: indeed, $\qs A\in R,\ \qs f\in \co(G)$, $f$ 
is not constant on $A$ because $A$ contains at least one edge, so the
restriction of $f$ to $A$ is a well-defined element of
$C^A$. Moreover, this is an imbedding because
$V\in R$.

For any $g\in\gr$, its stratum will be defined by the part $S$ of
$R$ described in one of the two following equivalent ways.

First, the elements of $S$ can be enumerated recursively by the following 
algorithm, starting
with its first element $A=V$:  
Consider the partition of $A$ defined by $g_A$
(i.e. the set of preimages of the singletons). Each of the 
connected components of these parts with cardinality greater
than one will be a new element of $S$, and the same operation 
must be done with this element.

Second, $S$ is the set of all
$A\in R$ such that for all $B\in R$ that strictly contains
$A$, 
$g_B$ is constant on $A$. 

Now we describe these strata $\ss$ and their imbeddings in $\h$.

\begin{nota}
$G/A$ denotes the graph obtained by identifying the elements of $A$ 
into one vertex%%denoted $\a $
, and deleting the edges in $E_A$;

Let $S$ be any part of $R$ such that $V\in S$ and any two elements 
of $S$ are either disjoint or included one into the other. 

For all $A\in S$, $A/S$ denotes the graph $A$ quotiented (like above) 
by the greatest elements $A_i$ of $S$ strictly included
in $A$ (they are disjoint and any element of $S$ strictly included in $A$
is contained in one of them).
Let$$\ss=\prod_{A\in S} \co(A/S).$$
Denote $S'=S-\{V\}$.

$\qs A\in R$, $\ba$ denotes the smallest element of $S$ containing $A$;

$\qs A\in S'$, $\ha$ denotes the smallest element of $S$ that strictly 
contains $A$
\end{nota}

Let us first check that the last notations make sense. 

For any nonempty part $A$ of $V$, the smallest element of $S$ containing $A$
is well-defined because first $A\inc V\in S$, then two elements of 
$S$ which contain
$B$ cannot be disjoint, so one of them is included in the other.

The definition of $\ha$ is justified in the same way. This gives $S$ a
tree structure.

The space $\co(A/S)$ will be identified with a part of $C^A$: the
$f\in\co(A/S)$ are then constant on each of the sets $A_i$.

Note that the interior stratum, $\co(G)$, is just $\co_{\{V\}}$.

\begin{lemma}\label{plo} There is a canonical imbedding of $\ss$ in 
$\h$ defined as follows: the image of an element
of $\ss$ in a $C^B$($B\in R$) is the restriction of its part
in $\co(\bb/S)$ to $B$. This way, the space $\gr$ is the disjoint
union of the $\ss$ where $S$ runs over the
parts $S$ of $R$ such that $V\in S$ and any two elements of $S$ are either 
disjoint or one included in the other. 
\end{lemma}

Let us prove this imbedding is well-defined (the second part of the lemma 
is easy).

For any $B\in R$ and any $f\in \co(\bb/S)$, we would like the
restriction of $f$ to $B$ to define an element of $C^B$. So we have to show
that $f$ is not constant on $B$.

Note that $B$ can  be included in none of the elements of $S$ which are
strictly included in
$\bb$, for $\bb$ is the smallest element of $ S$ containing $B$.
So $B$ is not reduced to a vertex in $\bb/S$, and since $B$ is connected,
one of its edges is an edge of $\bb/S$. This proves
that $f$ is not constant on $B$, by definition of $\co(\bb/S)$.
Finally, 
$\ss$ imbeds in $\h$ because $S$ is included in $R$.

{\it Note: }If $B\in R$, $\bb$ is also
the greatest element (and the union of elements) $A$ of $R$ containing $B$ 
such that $g_A$ is not constant on $B$.

\subsection{Description of the corners in $\gr$}

A direct calculation shows that dim 
$\ss=\dim\gr-\#S'$. We are going to see how $\ss$ has precisely
a neighbourhood\footnote{A neighbourhood of
a part  $P$ is a neighbourhood of every point  of $P$; 
it does not necessarily contain the closure of $P$.} 
in $\gr$ which is diffeomorphic to a neighbourhood of 
$\ss\times\{0\}^{S'}$ in $\ss\times[0, \1[^{S'}$. This diffeomorphism
is not completely canonical: to define it we have to choose an 
identification of $\co(A/S)$ 
with an arbitrary smooth section of representants in $(\ra)^{A/S}$, where
$\ra$ is just a copy of $\R^3$ marked with the label $A$, 
for each $A\in S$; we suppose that the origin of $\ra$ always
belongs to the affine space generated by the image of $A$.

Then, to a family $u=(u_A)_{A\in S'}$ of positive real numbers and a family
$f\in \ss$ (that is, $f=(f_A)_{A\in R}$ and $f_A\in \co(A/S)$), associate
the net of correspondences from each space $\ra$ for $A\in S'$ to the space
$\rb$ defined by
$$\eqalign{\phi_A:\ra&\longrightarrow\rb\cr x&\longmapsto f_\ha(A)+u_Ax\cr}$$

Then, define an element $g$ of $\h$ in the following way:

For all $B\in R$, to define the image of $g$ in $C^B$ (as a map 
from $B$ to $R_\bb$), first take for all $v\in B$ its image
by
$f_{\overline{\{v\}}}$ in $R_{\overline{\{v\}}}$, then compose all the
maps $\phi_A$ above for the $A\in S$ such that
$\overline{\{v\}}\inc A\mathrel{\raise 1pt
\vtop{\hbox{$\subset$}\kern -6.7pt\hbox
{$\,\scriptscriptstyle\not=$}
}}\bb$ 
to obtain an element of $R_\bb$.

This map is not constant for small enough values of $u$ because for
$u=0$ it is equal to $f_\bb$, which has been seen (in the proof
of Lemma \ref{plo}) to be not constant
on $B$.

This provides a well-defined smooth map from a 
neighbourhood $N$ of $\ss\times\{0\}^{S'}$ in $\ss\times[0, \1[^{S'}$, to a
subset of $\h$, and its restriction to $\ss\times\{0\}$ is the imbedding
defined in Lemma \ref{plo}.

It maps the interior of $N$ to $\co(G)$, because when $\qs A\in S'$, $u_A>0$,
these operations are equivalent to first taking the image of $g$ in $C^V$,
then restricting it to $A\in S$, because all the maps $\phi_A$ 
are dilations-translations of
$\R^3$. So it maps $N$ to $\gr$.

One easily checks in the same way that it maps any $(f, u)\in N$ to 
the stratum $\co_{S_u}$ where $S_u=\{A\in S|A=V\hbox{ or }u_A=0\})$.

Conversely, let us check that the $u_A$ ($A\in S'$) and the $f_A$ ($A\in S$)
can be expressed as smooth functions in a neighbourhood of $\ss$ in $\h$.
Let us first look at the case $S=\{V,A\}$.
So let $g$ be an element of $\gr$ in a neighbourhood of  $\ss$.
We shall more precisely use its components $g_V$ and $g_A$. First, we have 
$f_A=g_A$. Then, we define barycentric coordinates of the origin of $R_A$ 
with respect to the images by $g_A$ of the elements of $A$, 
as smooth functions of $g_A$. These coordinates, when applied to the family
of images by $g_V$ of the elements of $A$, give a point of
$R_V$. Now, $f_V$ will map $A$ to this single point, and will be equal
to $g_V$ on $V-A$.

Now we have to show that the so-defined $f_V$ is not constant: since
$g$ is near $\ss$, the diameter of $g_V(A)$ is much smaller than the
diameter of $g_V(V)$, so $g_V(V)$ contains a point far from $g_V(A)$ compared
to the diameter of $g_V(A)$, thus different from its barycenter.

Now we have defined $f_V$; it takes values in $R_V$, so we define
$u_A$ to be the rate of the unique translation-dilation $\phi_A$
from $R_A$ to $R_V$ which transports $g_A$ to the restriction of $g_V$ to $A$.

The general case for $S$ can be deduced from this case by induction.
\cqfd
\medskip

The description of $\gr$ is now finished and we can proceed with the 
compactification of $\c$. 

\subsection{The compactified $\cb$ of $\c$}

Starting with a diagram $\g\in\dl$, define the graph
$G$ with the same set $V$ of vertices as $\g$, and define the set 
$E$ of edges as being
the set of visible edges $E^v$, plus the set of all pairs of consecutive
univalent vertices $\{x,\s(x)\}$ for $x\in U$.

First, canonically imbed the space $\c$ into the space

$$\hh=\h\times(S^1)^U$$
and define the compactified $\cb$ of $\c$ to be its closure in $\hh$. 

Now we are going to describe the corners of this space. 
The stratum of an element 
$(f,f')$ of $\cb$ is classified by the data of the stratum $\co_S$ of its 
part $f$ in $\h$, together with the partition $P$ of $U$ defined by its part 
$f'$ in $(S^1)^U$. 

\begin{prop} A couple $(S, P)$ belongs to the image of $\cb$ by the above
map if and only if: 
% $(S^1)^U\times \co_S$ rencontre $\cb$ si et seulement si 
For all $A\in S$, $A\cap U$ is an
interval for the cycle $\s$, and $P$ is of one of the following two forms:

(i) the trivial partition into one part 

(ii) the partition defined by $\bu/S$. 
\end{prop}

Note that in case (ii), the vertices in $V-\bu$ are those which
escape to infinity, or for which any path to a univalent vertex
passes through trivalent vertices which escape to infinity 
(in case (i) there can be an undetermination).

Then, the elements $(f, f')$ of this stratum defined by $(S, P)$
are the elements of 
$\co_S\times(S^1)^U$ which verify the following conditions: 

First, $f'$ must of course give the partition $P$ again, and respect 
the cyclic order $\s$. 

Second, in case (ii), $f_\bu$ must coincide on $U$ with $K\circ f'$
(we identify $R_\bu$ to the space $\R^3$ in which the knot is imbedded).

Third, for the other cases of an $A\in S$ such that $\#(f_A(A\cap U))\geq2$ 
(in other words, when the image of $A\cap U$ in $A/S$ has at least two
elements but $f'(A\cap U)$ is a singleton), then $f_A$ maps $A\cap U$ 
to a straight line $\D$ directed and oriented by the tangent vector 
to the knot at $f'(A\cap U)$, preserving its order: we put on $A\cap U$
a linear order which is ``the restriction
of $\s$ to $A\cap U$ viewed as a subinterval of $U$''. There are 
$\#(U/S)$ possibilities for this linear order if $A=\bu$ (in case (i))
and only one otherwise. Here $U/S$ means the image of $U$ in $\bu/S$, 
and it has at least two elements. 

By convention (for the smooth sections by which the spaces $\ra$ are seen as
copies of $\R^3$), we shall suppose $\D$ contains the origin of $\ra$,
and that if $f_A(A\cap U)$
is a singleton, % (this can be translated to a condition on $S$), 
it is at the origin.

Now we describe the neighbourhood of this stratum. This is again
a corner, with a family of independent positive
parameters which are the same as before (indexed by $S'$), plus one 
more in the case (i), which will be denoted $u_0$. 

To be quick, we shall only define the elements of $\c$ which correspond to
a family of strictly positive coordinates $(u_A)_{A\in S'}$, and $u_0$. 

This will be made in two steps: first an approximate definition, then
a correction. 

In the first step, we start with the family $(u_A)$ and an element of 
$\ss$, and the same construction as above gives an identification of 
all the spaces $\ra$ together (because $u_A>0$ for all $A$), 
and thus an element of $\co$. There remains to
identify one of the spaces (namely $R_\bu$) to the space $\R^3$ in which 
the knot is imbedded. 
In case (ii), this identification is
already done, whereas in case (i), it will be $$x\mapsto K\circ f'(U) +u_0x.$$

For the second step, we have to choose for each $\al$ in the image of
$K\circ f'$ a diffeomorphism $\ph_\al$ between two neighbourhoods of $\al$ 
which verifies the two conditions: first, it
approximates the identity near $\al$ up to the first order; second, it maps
the tangent line to the knot at $\al$, to the knot itself. This system
of diffeomorphisms must depend smoothly on $f$.

Then for every $\al$, correct the map from $V$ to $\R^3$ of the first step 
by composing it with the map $\ph_\al$ for the vertices in 
${f_\bu}^{-1}(\al)$.

It can easily be seen that the resulting map extends as a smooth map
on a %smoothly to a 
neighbourhood of the stratum.

\section{List of the faces}
\setcounter{equation}{0}

Now we can list the codimension 1 strata of $\cb$. They are those
for which there is only one coordinate. Let us classify them into five types:

(a) In case (i) there is the coordinate $u_0$, so $S$ must be reduced to
$\{V\}$.
\smallskip

In case (ii), $S$ is of the form $\{V,A\}$, with four
possibilities:

(b) $U\inc A$ (this corresponds to the case when the vertices not in $A$ go
to infinity),

(c) $A\cap U=\0$,

(d) $A$ is of none of the types (b) and (c), and $E'_A\not=\0$,

(a') $A$ is of none of the types (b) and (c), and $E'_A=\0$.

\begin{nota} We denote the face of $\cc$ defined by $S=\{V,A\}$ by $\B (\g,A)$.
We say that a face is {\it degenerate} if its image by $\Phi$ is of
codimension strictly greater than one in $(\sd)^E$.
\end{nota}

The aim of this section will be to eliminate some types of faces by showing
they are degenerate, and obtain the following proposition:

\begin{prop}\label{nondeg} The only faces of $\cb$ which can be non-degenerate
are those of type (a) with $\g$ connected, and the faces $\B (\g,A)$ 
of the following forms: 

Those of type (a') with $A$ connected;

Those of type (c) or (d) such that $\#A=2$ (which will correspond to the 
IHX and STU
gluings);

Those of type (c) or (d) such that $\#A>2$, $A$ is connected and each 
trivalent vertex in $A$ meets at least two edges in $E_A$, and such that
$\#E'_A=1$
or $2$ for the type (d), $\#E'_A=3$ or $4$ for the type (c).
\end{prop}

These faces of types (a) and (a') are the ``anomaly" faces considered by Bott
and Taubes, for which we will have to introduce new spaces in Section 6.
The other faces will be glued in Section 5.

First, %etailing them, 
let us make a remark.

\begin{remark} 
Suppose there exists a part $A$ of $V$ such that $E'_A=\0$. Then
consider the two subdiagrams $\g_1$ and $\g_2$ whose sets of vertices are
$A$ and $V-A$ respectively. %, with the proper cycles.
Then there is a canonical diffeomorphism from $\c$ to an open part of
$C_{\g_1}\times C_{\g_2}$, which is 
delimited by a finite union of subspaces. 
\end{remark}

This remark has the two following consequences:

\begin{prop}\label{group}
$Z(K)$ is a grouplike element of $\A$.
\end{prop}

and

\begin{lemma}\label{tag} A type (a) face with $\g$ not connected, or a type (a') face
with $A$ not connected, is degenerate.
\end{lemma}

{\it Proof of Proposition \ref{group}} To prove that $Z(K)$ is grouplike, 
we first reformulate the expression
of $Z$. We define the notion of {\it knot graph} as Yang did.
In a knot graph, the set of vertices is defined to be a finite
part of $\rt$. Then an edge is a pair of vertices, the univalent vertices
are on the knot, and so on.

For a diagram $\g\in \ed$, the space $\cem$ of the knot graphs which
are isomorphic to $\g$ is the quotient space of our configuration space 
$\c$ by the automorphism group of $\g$. So we would like to consider the 
expression $$I_K(\g)[\g]\over|\g|$$ as an integral over the space $\cem$ of
the quotient form of $\Phi^*\O$. The problem is that this space $\cem$
may be non-orientable. But if we replace the differential form $[\g]\Phi^*\O$
by the corresponding measure on $\c$, then the quotient is well-defined.
Indeed, this measure does not depend on labellings or orientations,
so it is invariant under automorphisms.

Now, the proof that $Z$ is grouplike is straightforward.
The expression $Z\tens Z$ is the integral over the couples of knot graphs
$(\g_1,\g_2)$ with coefficient the tensor product 
$[\g_1]\tens [\g_2]$. This is the same as the integral over the knot graphs 
$\g$ with coefficient the sum of the 
$[\g_1]\tens [\g_2]$ for the couples $(\g_1,\g_2)$ such
that $\g=\g_1\sqcup\g_2$. This coefficient is just $\Delta([\g])$.\cqfd

\medskip

{\it Proof of Lemma \ref{tag}.}
The remark gives a diffeomorphism between a type (a') face, 
and the product of a type (a) face of $C_{\g_1}$ with
the space $C_{\g_2}$. So it shows that the types (a) and (a') are fundamentally
the same.

Now, for the type (a), if $\g$ is not connected, we see that
the fiber over an $f(U)\in S^1$ is mapped to the product of fibers
corresponding to the different subdiagrams, and this loses at least two
dimensions (which are the dilations and translations of one component
relatively to the other).

This implies also the result for the type (a').\cqfd

\begin{lemma}
Any type (b) face is degenerate (and the codimension of its image is greater
than 3).\end{lemma}

Set $E^c_A=E^v-E_A$. A type (b) face is the cartesian product 
of $\co(V/A)$ with the space $C_A$, and $\Phi_\g$ is the direct product
of the two maps
$$\co(V/A)\lra(\sd)^{E^c_A}$$and $$C_A\lra (\sd)^{E_A}.$$

Now, when applying (\ref{count}) to the complement of $A$ we find
$$\dim\co(V/A)=\dim(\sd)^{E^c_A}-\#E'_A-1.$$
Finally, since $\g$ is triply connected to $U$ we have $\#E'_A\geq 3$.\cqfd

\begin{lemma}
A type (c) face can be non-degenerate only if 
$\#E'_A=3$ or $4$; a type (d) face can be non-degenerate only if $\#E'_A\leq2$.
\end{lemma}

In case (c), $\Phi$ is the product of the two maps:
$$C_{\g/A}\lra(\sd)^{E^c_A}$$ and $$\co(A)\lra (\sd)^{E_A}$$
But $$\eqalign{\dim(\sd)^{E^c_A}&=\dim C_{\g/A}+\#E'_A-3\cr
\dim(\sd)^{E_A}&=\dim \co(A)+4-\#E'_A.}$$

For the type (d) we consider the map $C_{\g/A}\lra(\sd)^{E^c_A}$,
for which $$\dim(\sd)^{E^c_A}=\dim C_{\g/A}+\#E'_A-1.$$
\cqfd

\begin{lemma}
A type (d) face can be non-degenerate only if 
$A$ is connected or $\#A=2$.
\end{lemma}

Indeed, if $\#A>2$ and $A$ is not connected, then for a fixed direction $\D$,
the map from $C^A$ to
$(S^2)^{E_A}$ factorises through the product of the $C^B$ where $B$ runs
through the connected components of $A$, so we lose at least one dimension
in the relative translations and dilations of the components.
\cqfd

\begin{lemma}
Any type (c) or (d) face such that $\#A>2$, $\exists B\inc A\cap T,\ 
\#(E'_B\cap E_A)=1$, is degenerate.
\end{lemma}

We suppose that $A$ is connected, thanks to the preceding lemma.
Then there exist two parts $A_1$, $A_2$ of $A$
with $A=A_1\cup A_2$, $\#(A_1\cap A_2)=1$, $\#A_1>1$, $\#A_2>1$, 
$A_2\inc T$ and
$E_A=E_{A_1}\cup E_{A_2}$. If $\#B>1$, we take $A_2=B$.
Else, we take for $A_2$ the two ends of the edge in $E'_B\cap E_A$.
The map from $C^A$ to $(S^2)^{E_A}$ factorises through
$C^{A_1}\times C^{A_2}$. This loses one dimension (the relative
scale between the two configurations). What we use here
is the following argument.
The map $\Phi$ maps the face  $\B (\g,A)$ with
dimension $2N-1$ into the space $(\sd)^E$  with dimension $2N$. But the one-parameter
group of dilations of $A_2$ relatively to 
$A_1$ acts freely on $\B (\g,A)$, and the function
$\Phi$ is constant on its orbits, which proves that
$\dim \Phi(\B (\g,A))<\dim 
\B (\g,A)=2N-1$.\cqfd
\medskip

We deduce that if a type (c) or (d) face such that $\#A>2$ is non-degenerate, 
then any trivalent vertex in $A$ has at least two edges in $E_A$.
This ends the proof of Proposition 4.2.

\section{Ordinary gluings}\label{recol}
\setcounter{equation}{0}

The aim of this section is to prove the following proposition:

\begin{prop} Let $\a $ be the map from $\dl$ to $\A_n^k$ defined by
Formula $(\ref{coef})$:
$$\a (\g)={(u_\g-k)!\over (3n-k)!\;2^{3n-u_\g}}[\g].$$
There exists a gluing of the faces of the linear combination of spaces
$$C=\sum_{\g\in \dl}\a (\g) \cc,$$ compatible with $\Phi$, such that the only
remaining non-degenerate faces are the anomaly ones.\end{prop}

This is the consequence of the three following lemmas, which will be proved
in Subsection 5.2. They concern a map $\a $ from $\dl$ to any 
$\Z$-module, which is extended to the diagrams in $\dl$ with any other
orientation with respect to the relation AS:

\begin{lemma} The gluings of the type (c) faces with $\#A=2$ are ensured
by the following RIHX relation between labelled diagrams that are 
identical outside the drawn part, and all belong to $\dl$:
$$\a \epfb{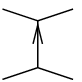}+\a \epfb{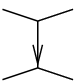}=\a \epfb{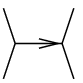}+\a \epfb{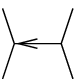}-\a \epfb{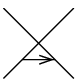}-\a \epfb{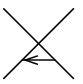}$$
\end{lemma}

\begin{lemma} The gluings of the type (d) faces with $\#A=2$ are ensured
by the following RSTU relation between labelled diagrams that are 
identical outside the drawn part:
$$\a \epfb{u}-\a \epfb{s}
=\sum_{e\in E^a}\left( \a \epfb{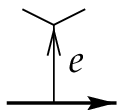}+\a \epfb{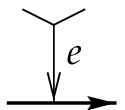}\right).$$
where $E^a$ is the set of absent edges of a diagram in the left-hand side 
of the relation. This relation occurs when all diagrams are
triply connected to $U$, and the left-hand side diagrams belong to $\dl$.
(the right-hand sum can be empty, but the unlabelled diagrams it would
involve must be triply connected to $U$).
\end{lemma}

\begin{lemma} The type (c) or (d) faces with $\#A>2$ 
are glued when $\a $ is a constant on each isomorphism class of diagrams.
\end{lemma}

\begin{remark} The above conditions are sufficient conditions for the gluings,
and not necessary conditions. To prove Theorem \ref{mult}, we just
need to prove that the formula (\ref{coef}) satisfies these conditions.
We shall see in Subsections 7.3 and 7.4 that other choices of 
formulas for $\a $ can provide better
rationality results that the ones of Theorem \ref{mult}.
\end{remark}

\subsection{Gluing methods and orientations}

We are going to define gluings of faces, that is, bijections $b$ between
two faces $\B (\g,A)$ and $\B (\g',A')$ such that $\Phi\circ b=\Phi$ on $\B (\g,A)$. 
They generate an equivalence relation between faces, where two faces are
equivalent iff there exists a gluing between them.

Then, the (algebraic) gluing condition on an equivalence class of
faces $\B (\g,A)$ is that$$\sum_{(\g,A)} \pm \a (\g)=0.$$
Here the signs of the terms are determined by the following rules.

As a face $\B (\g,A)$ inherits an orientation from the orientation 
of the space $\c$, the signs in front of $\a (\g)$ and
$\a (\g')$ for two given couples $(\g,A)$ and $(\g',A')$ will be equal iff
the identification between the faces $\B (\g,A)$ and $\B (\g',A')$ preserves
their orientations. (If there is no distribution of signs obeying this
rule for all gluings inside an equivalence class, then each face is glued to
itself, so there is no condition anymore).

Now we are going to define two general cases of gluings between two faces
$\B (\g,A)$ and $\B (\g',A')$ of a same type ((c) or (d)). The problem will be
to compare the orientations of the two faces. Then, these
constructions will be applied to several more particular cases.

We suppose that $\g/A=\g'/A'$ (as labelled diagrams), and that the gluing
corresponds to the identity on $\co(V/A)$. 
So, the only transformations concern the graphs $A$ and $A'$ and the 
configuration spaces $\co(A)$ and $\co(A')$, and, for the (d) type,
the univalent vertices of $A$ and $A'$ are mapped to the same line $\D$ 
in the two identified configurations.

The two (but not disjoint) general cases that we shall consider will be
the following:

1) $E_A=E_{A'}$ (this implies that $A$ and $A'$ have the same respective
numbers of univalent and trivalent vertices).

2) $\#A=2$ ($=\#A'$, for $\g$ and $\g'$ have the same degree).
\medskip

{\bf Let us see the first general case.} So, $A$ and $A'$ have 
the same respective
numbers of univalent and trivalent vertices.

For convenience, we are going to label all vertices in $A$, starting with
the univalent vertices (this means the elements of $U$), if they exist.
We denote them by $x_i$ for $i\in J$, $J=\{0,\cdots\#A-1\}$. By convention,
we put the vertex $x_0$ at the origin of the space $R_A$ (it is
univalent in case (d)). So $\co(A)$ is identified with the space of half
lines in the vector space $(\R^3)^{\#A-1}$. 
We do the same for $\g'$ (so the univalent vertices in $A$ and $A'$
are labelled by the same elements of $J$). Precisely, we shall
suppose (for the (d) type) that over each $\D$, the bijection between 
the faces is defined
by a linear transformation of the subspace $\cal E$ of $(\R^3)^{\#A-1}$ 
consisting of the configurations which map $A\cap U$ to $\D$.

The neighbourhoods of the faces are locally identified with the cartesian
product $\ce\times C_{V/A}$. Then, the identification between them
reverses the orientation of $\ce\times C_{V/A}$ iff det $L>0$.

Now, let us deduce its effect on the orientations of $\c$ and $\c'$ as
defined in Section \ref{ori}: the local coordinates, first labelled by
the vertices, are then labelled by the half-edges. So this labelling
is changed when we change the structure of a graph, or just its labelling.

So, the graphs $A$ and $A'$ define two bijections $b$ and $b'$ 
between
$$\{i\in J|x_i\in U\}\cup\{(i,1),(i,2),(i,3)|i\in J,x_i\in T\}$$
and the set $\cup A\inc X$ of half-edges which belong to a vertex in $A$.
Then, the sign between the orientations of $\B (\g,A)$ and $\B (\g',A')$ is
the product of the sign  of det $L$ with the signature of $b^{-1}b'$.

{\bf Let us see the second general case.} 

Now $\#A=2$, and the face $\B (\g, A)$ can be 
identified with the space $C'_{\g/A}$. 
Indeed, let us distinguish two cases.
If $E_A=\0$, then $\B (\g, A)$ is identified with $C_{\g/A}$ (since the
configuration of $A$ is determined by the one of $\g/A$, according to the
tangent map of the knot). If $E_A\not=\0$, we have
$A\not\inc U$ (since we are only interested in types (c) and (d) faces).
So $\B (\g, A)$ is the product of $C_{\g/A}$ by $S^2$. This allows us to put
the unique element of $E_A$ among the absent edges (since $\Phi$ is
factorised by the identity on this $\sd$, with well-chosen conventions).

\begin{remark} The image of $A $ in  $\g/A$ has one more half-edge 
than an ordinary vertex:
four if it runs through $\rt$, and two if it runs along the knot. 
The additional
half-edge represents the inward normal (i.e. the parameter $u_A$).
This defines an orientation on the normal to the image of a non-degenerate
face $\B (\g,A)$ in $(\sd)^E$. 
\end{remark}

\subsection{Proof of Lemmas 5.2 and 5.3}

Note that in all cases, all diagrams involved in an RIHX or RSTU relation 
are distinct as labelled diagrams.

Let us first see the case when all the diagrams 
involved in an RIHX or RSTU relation 
are represented in $\dl$.

The RIHX relation is both an application of the first and second general
cases. We deduce it from the following remark.

In the first general case, we can always take $L=$ Id, preserve the
graph structure of $A$ but just use permutations of the set $E'_A$ 
(modulo the transpositions of two edges attached to the same vertex
of $A$). This provides non-trivial gluings
when the edges in $E'_A$ are attached to different vertices in $A$.

The RSTU relation is an application of the second general case. To check the
signs, let us study the gluing between two faces corresponding to diagrams
with a different number of univalent vertices (as the other ones can
easily be deduced from the first general case).

Consider the following diagrams, where each half-edge is marked with the
name of a component in a local coordinate system of the neighbourhood
of each face.
Here we use the spherical coordinates for the trivalent vertex:
$$\epsfbox{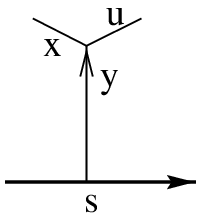}\hskip 3cm\epsfbox{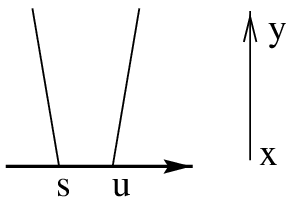}$$

We see that these two labellings differ by the transposition $(x,s)$.
This proves that these two diagrams appear in the relation with different
signs (or on different sides of the equality).

\medskip

{\bf Study of the exceptional cases}

Now, let us face the possibility that some diagram in an  IHX or STU relation
is not represented in $\dl$. There are two situations.

First, suppose one of the diagrams is not triply connected to $U$.
Such a relation must be simply deleted because these faces are degenerate,
since one of them is (according to Lemma 2.10), and $\Phi$ coincides on them.

Second, suppose all diagrams in the STU relation are triply 
connected to $U$, but not all are represented in $\dl$. This implies
that two diagrams have $k$ univalent vertices and the other one $\g$
has $k-1$ ones.
This implies for the corresponding RSTU relation that the sum on the 
right-hand side is empty. So there is no problem: this sum is zero,
and this relation is verified for the expression (\ref{coef}) because
it is deduced from the STU relation and the cancellation of 
$[\g]$ in $\A_n^k$.

This ends the proof of Lemmas 5.2 and 5.3. \cqfd

\subsection{Proof of Lemma 5.4}

We consider now a non-degenerate type (c) or (d) face with $\#A>2$.
So, we are going to apply the first general case.

According to Proposition \ref{nondeg}, any vertex $x\in A$ attached to
an arbitrary edge in $E'_A$ is trivalent (because $A$ is connected),
and its two other edges $e_1$ and $e_2$ are in $E_A$. Let $y$ and $z$ be 
the other
ends of these edges. Give $x$ a non-zero label, and let $L$ be the map
which applies to $x$ the central symmetry with respect to the
middle of $y$ and $z$, that is:$$x\longmapsto y+z-x,$$and preserves the
other vertices.
This transformation has determinant $-1$.

Then apply the only change of structure of $A$ which fits with this 
transformation: it is a mere relabelling, which exchanges $e_1$ and $e_2$ 
and reverses them: if $e_1$ goes from $x$ to $y$ in $A$, then it
goes from $z$ to $x$ in $A'$. 
Suppose these diagrams have the same orientation through this relabelling.
This defines an even permutation of $\cup A$, so the two faces
have an opposite orientation in the gluing.

Now, we just have to prove that these faces can be glued pairwise this way
in the particular case when $\a$ is constant on each isomorphism class of
diagrams with the same orientation.

Consider an (equivalence) class of couples $(\g,A)$ (concerned here) 
such that $\g/A$ is fixed as a labelled diagram, and the
isomorphism class of the couple $(\g,A)$ is also fixed. Then fix 
an arbitrary edge $e\in E'_A$. 
By this choice, the diagrams in this class are grouped in pairs
according to the above construction.\cqfd

\section{Anomaly}
\setcounter{equation}{0}

\subsection{The spaces $W(\g)$}

Remember that the non-degenerate types (a) and (a') faces are defined by
a connected component $A $ of $\g$, whose configurations map all
univalent vertices of $A$ to a line $\D$ directed 
by the tangent vector $x$ to the
knot. So let us give a few definitions:

\begin{defin} Let $\dw$ 
be the set of diagrams 
$\gu=(U,\leq,T,E^v)$ defined as in Definition \ref{diag}, 
except that: the cyclic order $\s$ over $U$ is replaced by a linear
order $\leq$, and we suppose that $\gu$ 
is connected.
Once more, each of these diagrams has an arbitrary orientation.

For all $x\in S^2$ and $\gu\in\dw
$, $\wx(\gu)$ denotes the space of maps $f$ from the set $V$ 
(of vertices of $\gu$) to $\rt$, such that: f is non-constant on each edge, 
and maps $U$ to the line
$\D$ oriented by $x$, preserving the order of $U$. And $f$ is considered
modulo the translations and dilations which preserve $\D$.

Then, $W(\gu)$ denotes the union of the $\wx(\gu)$ for $x\in S^2$.
\end{defin}

We are going to give these spaces an orientation. First, orient
the two-dimensional Lie algebra of the group of translations 
and dilations which preserve $\D$, by giving it a base: take 
a first vector oriented to the positive translations of the oriented axis
$\D$, and a second vector oriented to the dilations of rate greater than one.

Now, since the space of maps $f$ from $V$ to $\rt$ which map $U$ to
$\D$ is oriented by the method of Subsection \ref{ori}, 
this orients the quotient space $\wx(\gu)$. Moreover, since this space is
a fiber for the fibered space $W(\gu)$ with basis $\sd$, the space
$W(\gu)$ is also oriented. There is no ambiguity, since all
considered spaces are even-dimensional.

Then we define a map  $\Phi_\gu$ from $W(\gu)$ into 
$(\sd)^{E^v}$, and a map $\Phi'_\gu$ from 
$W'(\gu)=W(\gu)\times (\sd)^{E^a}$ into $(\sd)^{E}$, exactly as before.

Now we are going to see that these spaces can be glued together according to
the same formulas, except that there are no remaining faces
such as the types (a) and (a') faces.

\begin{prop}\label{w} Suppose a map $\a $ from $\dw$  
to a
$\Z$-module verifies any relation of Lemmas 5.2, 5.3, 5.4, when
all the diagrams in the relation are connected.
Then the linear combination of spaces
$$\sum_{\gu\in\dw
}\a (\gu) W'(\gu)$$ is completely glued, hence the union $\Phi$ of the
$\Phi'_\gu$ admits a degree on it.
\end{prop}

{\it Proof}

The method is to use the same proofs as those of Sections 4 and 5 in 
this new situation, where there are no equivalents of types (a) and (a') faces
(for the scale is relative and the diagrams are connected)

We just have to check that the map $\Phi$ defined on the spaces
$W(\gu)$ factorises the same way, and can be decomposed the same way
on the faces, the dimensions being the same.

It is very easy to do for the faces $\B (\gu,A)$ such that $A\inc T$:
$$\B (\gu,A)=W(\gu/A)\times \co(A).$$

The other faces can also be seen as fibered spaces with basis
$W(\gu/A)$ and fiber $W_x(A)$, for the value of $x$ determined by the
element of $W(\gu/A)$.

The only new thing to prove is that if a relation involves a non-connected
diagram, then the faces are degenerate. This can easily be done
by splitting the fiber over $\sd$ into the product of the corresponding 
fibers for the connected components. This loses two dimensions,
corresponding to the relative translations and dilations between the
configurations of the two subdiagrams.
\cqfd

\begin{defin} Let $\ewn$ be the set of isomorphism classes of degree $n$
connected diagrams, with a linear order on univalent vertices, and with
arbitrary orientations on the trivalent vertices. 
For all $\gu\in\ewn$ we define the 
integral
$$f_\gu=\int_{W(\gu)}\Phi_\gu^*\O$$

Altschuler and Freidel \cite{af} proved that $f_\gu=0$ when deg $\gu$ 
is even, for symmetry reasons.

Then we define the degree $n$ anomaly to be 
$$\an =\sum_{\gu\in\ewn}{f_\gu\over|\gu|}
[\gu]\in \A_n$$
where $|\gu|$ is the number of automorphisms of $\gu$ which preserve the linear
order on $U$, and where $[\gu]$ denotes the image in $\A$ of the diagram
obtained from $\gu$ by closing its linear order into a cyclic order.
Finally,$$\ano =\sum_{n>0}\an .$$
\end{defin}

Note that $\an $ is the degree of the map $\Phi$ of the 
above proposition, for the map $\a $ defined by the same formula
as (\ref{coef}) when $\A_n^k=\A_n$:
$$\a (\gu)={(u_\gu-k)!\over (3n-k)!\;2^{3n-u_\gu}}[\gu].$$
Moreover, $a_1=[\th]$.% (see the proof of Lemma \ref{tet})

\subsection{The generalised diagrams and their configuration spaces}
\label{gen}

Now we are going to define the additional configuration spaces which will allow
us to glue all the faces and define the degree of $\Phi$.

First, define a smooth map $\ph:\dd\longrightarrow
S^2$ from the disc to the sphere, which coincides with the tangent map 
d$K$ to the knot on the boundary $S^1$ of $\dd$.
Here we give the disc $\dd$ the opposite orientation to 
the usual one. So, the basis (tangent oriented by $S^1$, outward
normal) is direct.

\begin{note} For all $\gu\in\dw$ and $\ph:\dd\longrightarrow
S^2$, let $$\wg= \ph^*W(\gu)$$ be the fibered space over $\dd$ which is
the union of the
$W_{\ph(x)}(\gu)$ for $x\in \dd$, and put the product orientation on it.
\end{note}

Now we can define the generalised diagrams and their configuration spaces.

\begin{defin} A {\it generalised diagram} is defined as in (\ref{diag}),
where $\s$ is replaced by the following data:

\b A set $Y$ of connected components of $\g$, together with a linear 
order on the set of univalent vertices of each of them. 

For convenience, we shall use a notation with indices, with $Y$ as
the set of indices, although this set depends on $\g$ (and has a variable
number of elements). Let $U_\gu$ be the set of univalent vertices of $\gu$
for $\gu\in Y$, $\g_0$ be the remaining part of the diagram (it can be
empty), and let $U_0$ be the set of its univalent vertices.

\b A cyclic order $\s$ on $U_0$.

We denote by $\dg$ the set of these generalised diagrams of degree $n$.
An arbitrary orientation is chosen on each of them.

\end{defin}

\begin{defin} Let the configuration space of a generalised diagram be
$$\c=C_{\g_0}\times\prod_{\gu\in Y}\wg.$$
with the product orientation. This is again completed with the help of the
absent edges:$$\cc=\c\times (\sd)^{E^a}.$$
\end{defin}

The product orientation is well-defined because all spaces are 
even-dimensional. It is reversed when we change the orientation
of a single vertex.

Now we can glue all these spaces together:

\begin{prop}\label{tot} Let $\a $ be a map from $\dg$ into a 
$\Z$-module 
which verifies the following conditions:

\b The restriction of $\a $ to diagrams such that $Y=\0$ satisfies the
conditions of Lemmas 5.2,5.3,5.4.

\b For all generalised diagrams $\g$ and $\g'$ which have the same edges
and vertices and such that $\g'$ is obtained from $\g$ by inserting one
of the ordered $U_\gu$ inside the cycle $\s$, we have $\a (\g)=\a (\g')$.

Then the canonical map $\Phi$ defined on the linear combination
$$\sum_{\g\in\dg}\a (\g) \cc$$ admits a degree.

\end{prop}

{\it Proof}: 

{\bf 1) Ordinary gluings}

The ordinary gluings of $\g_0$ are obvious. As for the other
ordinary gluings (which involve faces of the spaces $W'(\gu)$), the two
above conditions on $\a $ imply the conditions of Proposition 
\ref{w}, where we fix $\g_0$, and all the $\gu$ except one which
varies freely. So here $\a $ is considered as a function of this last
diagram, and the RSTU relations involve all the absent edges of $\g$.
\medskip

{\bf 2) Anomaly gluings}

Now we are going to glue together the following types of faces.

\b The face of a $\cc$ defined by the boundary of the copy of $\dd$ labelled
by some $\gu\in Y$. This face is glued if $\g_0\not=\0$ to the following 
faces:

\b The type (a') faces $\B (\g',A)$ where $\g'$ runs over the diagrams obtained 
from $\g$ by inserting $U_\gu$ somewhere in the cycle $\s$, and $A$ is the
set of vertices of $\gu$.

\b But if $\g_0=\0$, it is glued to a unique type (a) face in an obvious way.

Now we have to check that the two orientations are opposite in any case.
The orientation of $\dd$ has been defined by the basis (tangent oriented
by $S^1$, outward normal), whereas for the type
(a) or (a') faces, we must consider the orientation of the Lie algebra
by which we quotiented, when we defined the orientation of the spaces
$W_x(\gu)$: we took the basis (translations in the sense of $x$, growing
dilations). But $x$ represents the generic tangent vector to the knot
with the positive orientation, whereas the growing dilations point inward
to the configuration space.
\cqfd\medskip

In particular, we shall consider the unique map $\a $ which extends
Formula (\ref{coef}) in this way:
$$\a (\g)={(u_\g-k)!\over (3n-k)!\;2^{3n-u_\g}}[\g]$$
The formula is the same, we just need to define the class of a generalised
diagram:
$$[\g]=[\g_0]\prod_{i\in Y}[\gu]$$
(Here we use the fact that the product is well-defined 
and commutative in $\A$.)

\subsection{The case when $I_K(\th)$ is an odd integer}

\begin{lemma} Define the integral $$I_\ph=\int_{\dd}\ph^*\om.$$
Then we have the identity$$\int_{\wg}\Phi_\gu^*\O_{E(\gu)}=f_\gu I_\ph.$$
\end{lemma}

This results from the invariance of the space $W(\g)$ by rotations of $\sd$.
\cqfd

Thus, the integrals of the configuration spaces of the new diagrams
(i.e. when $Y\not=\0$) will vanish provided that $I_\ph=0$. So, to
prove Theorem \ref{mult}, we just have to prove the following lemma:

\begin{lemma}\label{tet} 
If $I_K(\th)$ is an odd integer, then there exists 
$\ph$ such that $I_\ph=0$.
\end{lemma}

Let us study the configuration spaces of the generalised diagrams in degree
1. Take $E=\{e_1\}$. We have three generalised diagrams:
$\th$, $\th_+$ and $\th_-$. The two last ones correspond to the two possible
linear orders among the two univalent vertices. They are respectively
($\{1\}<\{2\}$)
and ($\{1\}>\{2\}$). These two diagrams are isomorphic.
We find that $f_{\th_+}=f_{\th_-}=1$ (we can check in the case of $\th_-$
that $\Phi$ is the antipodal map, but the fibre $W_x(\th_-)$, with dimension
0, has the negative orientation).

Then, consider the map $\a $ with values in $\Z$ defined by 
$\a (\th)=\a (\th_+)=\a (\th_-)=1$. According to Proposition 6.7, the
expression $$I_K(\th)+2I_\ph$$
is a constant integer on each isotopy class of the knot, when $\ph$ varies
continuously.

Now, let us specify the set of its possible values.
Start with the case of the unknot embedded as a circle in a plane. 
Then we have $I_K(\th)=0$ because $\Phi(C_\th)$ 
has dimension 1. Then, the possible values of $2I_\ph$ are the
odd (relative) integers.

Then, transform the unknot to any knot by a generic homotopy, that is, a
smooth map from $S^1\times [0,1]$ to $\rt$ such that for all
$t\in [0,1]$ we have an immersion from $S^1$ into $\rt$; this immersion
must be
an imbedding except for a finite series of values of $t$, where there is
precisely one self-intersection with two different tangent directions.
At these values, $\ph$ varies continuously but 
$I_K(\th)$ jumps two units. So, for any knot, the set of possible values of
$I_K(\th)+2I_\ph$ is always the set of odd integers. So, if
$I_K(\th)$ is an odd integer, then $0$ is a possible value for 
$I_\ph$.
This ends the proof of Theorem \ref{mult} in the case when $I_K(\th)$
is an odd integer.\cqfd

\subsection{The case when $I_K(\th)$ is an even integer}

\begin{prop}
If $I_K(\th)$ is an even integer, 
then the conclusion of Theorem \ref{mult} is also true.
\end{prop}

To prove this proposition, we shall express $Z_n$ as the degree of
a map defined on another linear combination of spaces.

Let $n'\leq n$ be an integer. Let $\g\in \dg$ and
$\gu\in Y$ such that $\deg\gu=n'$.

If $a_{n'}=0$, we keep the same construction, because the sum of these
terms is zero. Indeed, let us fix a labelled diagram except a $\gu$ of
degree $n'$ which varies freely (with a variable set of absent edges).
The map $\a $, as a function of $\gu$, 
obeys Formula $\ref{coef}$ (for another value of $k$), up to a constant
factor.

If $a_{n'}\not=0$ (in particular when $n'=1$), we are going to replace 
the spaces $\wg$ for $\deg\gu=n'$ by a twice lower quantity of spaces 
but which have
the same quantity of faces, with the help of the following constuction.
It will also prove again that each $f_\gu$ vanishes when deg $\gu$ is even.

\begin{lemma}\label{esp} For any integer $p$ there exists a compact oriented
surface $\Si$ with a two-component boundary, together with a smooth 
map
$\ph'$ from $\Si$ to $\sd$ and an involution $i$ of $\Si$ such that:

\b $i$ reverses the orientation of $\Si$

\b $i$ exchanges the two components of the boundary

\b $\ph'\circ i=-\ph'$

\b The restriction of $\ph'$ to one of the components of the boundary
coincides with the tangent map of $K$, and the orientation induced
by the orientation of $\Si$ of this component 
is the opposite of the orientation of $K$
(thus, its restriction to the other component is the opposite of the tangent 
map, with the orientation of $K$).

\b $I_{\ph'}+I_K(\th)=p$.
\end{lemma}

Let us first assume the lemma to prove the proposition.

We define the space $\wp$ in the same way as the space $\wg$.
If the diagrams $\gu$ and $\gu'$ are of degree $n'$, and if $\gu'$ is obtained
from $\gu$ by reversing the orientations of all the edges,
then we identify canonically the spaces $\wp$ and ${\cal W}_{\ph'}(\gu')$
through $i$ and the central symmetry of the whole configuration.
This identification is compatible with $\Phi$.

Let us see how this identification acts on the orientations of these
spaces.

We have a $-1$ sign for each trivalent vertex. There is no sign for each
univalent vertex, since the geometrical symmetry makes up for
the orientation change of the axis $\D$. And the central symmetry
of the configuration carries the orientation of the Lie group of 
translations and dilations by which we quotiented to obtain the space 
$\wx(\gu)$, into the orientation of this Lie group for $W_{-x}(\gu')$.
Moreover, the reversal of each edge induces a sign change,
and we have one more sign from the fact that $i$ reverses the orientation
of $\Si$. 

Finally, this identification between the spaces 
$\wp$ and ${\cal W}_{\ph'}(\gu')$ preserves the orientation if and only
if the degree $n'$ of $\gu$ is odd.

This first proves that the anomaly is zero at even degrees. Indeed, it
proves that if deg $\gu$ is even and if $\gu'$ is obtained from $\gu$ by
reversing all the edges, then $f_\gu=-f_{\gu'}$. But these diagrams are
isomorphic, so $f_\gu=f_{\gu'}$. Hence $f_\gu=0$.

Then, when $\anp \not=0$, $n'$ is odd and the above identification preserves
the orientation.

To avoid the repetition of these spaces, we suppose that for all generalised
diagrams $\g$ and $\g'$ which differ by reversing the orientations of
the edges on some
$U_\gu$ for $\gu\in Y$ and deg $\gu=n'$, we have $\a (\g)=\a (\g')$.
This is the case for Formula (\ref{coef}).
This relation between $\g$ and $\g'$ generates an equivalence
relation, and we are concerned with the quotient set of the set $\dg$ of
generalised diagrams, by this relation. 

Finally, if $I_K(\th)$ is an integer, put $p=-I_K(\th)$ to have 
$I_{\ph'}=0$. Now, to end the proof of
the proposition we only have to check the lemma:

\medskip

{\it Proof of Lemma \ref{esp}}

First case: $p=0$. Then, just take the surface $-C_\th$, 
with $\ph'=\Phi_\th$. The involution $i$ consists of exchanging
the positions of the two univalent vertices on the knot.

General case: take the disjoint union of the above solution with $p$
copies of $\sd$, with $\ph'=$~id and
$i=-$~id. \cqfd
\medskip
This ends the proof of Theorem \ref{mult}.

\section{Other results}
\setcounter{equation}{0}

\subsection{Explicit variation of $Z$}

We are going to make a new proof for the expression of the variations of
$Z$ on each isotopy class of knots, as a function of 
$I_K(\th)$, or equivalently of $I_\ph$.
(This formula had been proved by Altschuler and Freidel.)

\begin{prop} 
Let $n\geq 2$ and $k\leq 3$, so that $A_n^k=\A_n$. Then the expression 
$\deg \Phi$ of Proposition \ref{tot} applied to Formula \ref{coef}
is the degree $n$ part of $$Z \exp\left(I_\ph \ano\right).$$ 
\end{prop}

The proof of this proposition is just a matter of calculus. We have
$$\deg\Phi=\sum_{\g\in\dg}{(u_\g-k)!\over (3n-k)!\;2^{3n-u_\g}}
I(\g_0)[\g_0]\prod_{\gu\in Y}f_{\gu} I_\ph[\gu].$$

Let $\egn$ be the set of isomorphism classes of generalised diagrams of
degree $n$, and let $\ew$ be the union of the $\ewn$.
The number of automorphisms of a generalised diagram $\g$ which has 
$\beta_{\gu}^\g$ copies of $\gu$ for all $\gu\in \ew$ is
$$|\g|=|\g_0|\prod_{\gu\in\ew}|\gu|^{\beta_{\gu}^\g}\beta_{\gu}^\g!$$
\def\bb{\beta_{\gu}^\g}

So we have
$$\eqalign{\deg\Phi&=\sum_{\g\in\egn}
{I(\g_0)[\g_0]\over|\g|}\prod_{\gu\in Y}f_{\gu} I_\ph[\gu]\cr
&=\sum_{\g\in\egn}{I(\g_0)[\g_0]\over|\g_0|}\prod_{\gu\in\ew}
{1\over\beta_{\gu}^\g!}
\left({f_{\gu} I_\ph[\gu]\over|\gu|}\right)^{\beta_{\gu}^\g}\cr}
$$

This is the degree $n$ part of
$$\left(\sum_{\g_0\in\ed}{I(\g_0)[\g_0]\over|\g_0|}\right)
\exp\left(\sum_{\gu\in\ew}{f_{\gu} I_\ph[\gu]\over|\gu|}\right)$$
which is just 
$$Z \exp\left(I_\ph \ano\right).$$\cqfd

One can prove in the same way that, with the construction of Subsection 6.4,
$\deg \Phi$ is the degree $n$ part of $Z \exp({1\over2}I_{\ph'} \ano)$.

\begin{cor} On each isotopy class of knots, $Z$ is a function of $I_K(\th)$ of
the form $$Z=Z_0\exp\left({I_K(\th)\over2}\ano\right)$$
where $Z_0$ is an invariant of the knot $K$, and is rational at each degree.
\end{cor}

Indeed, take $p=0$ in Lemma \ref{esp}. This implies $I_{\ph'}=-I_K(\th)$. 
\cqfd

\subsection{Universality of $Z_n$ in the case of a free module}

We still fix an integer
$k\leq 2n$, and suppose $I_K(\th)$ is an odd integer. 
Recall that to any map $\a $ from $\dl$ to a $\R$-vector
space $\ce$ which verifies sufficient gluing conditions, we associate
a degree of the map $\Phi$ on the linear combination of the
configuration spaces:
$$\deg\Phi=\sum_{\g\in \dl}\a (\g) I_K(\g).$$

It would be very difficult to explicit necessary and sufficient gluing
conditions, because this problem has two aspects:

First, a topological aspect, which consists of making the complete list
of the degenerate faces, the possible gluings between the other faces, and
the configuration spaces which have a codimension 1 image.
We could also search for other spaces to add to this combination.

Second, a combinatorial aspect, where we make the equivalence classes
of faces explicit, and look for solutions of the resulting system 
which give better rationality results.

We could also wonder if some torsion Vassiliev invariant can be found by this
method, when the map $\a $ takes values in a torsion $\Z$-module.
But a torsion invariant, if it exists, might not be found by this method.
(This is suggested by the fact that the rationality results we can obtain in
degree 2 by this method are not the best ones: see Section 7.4). But
we won't be interessed in this problem here.

\begin{prop}\label{cond} The following system of sufficient gluing 
conditions can easily be deduced from the previous constructions
(when $I_K(\th)$ is odd):

1) The RSTU relation (see Lemma 5.3).

2) The RIHX relation (see Lemma 5.2), that we can replace by the 
following weaker RIHX' relation: a RIHX' relation is the sum of the
($\#E^a+1$) RIHX relations for which $\g/e$ is a fixed labelled diagram,
where $e$ is the edge involved.
This RIHX' relation is the one obtained as an application of the second
general case of Section 5.1.

3) The (complicated !) relations between the coefficients in each
isomorphism class of diagrams, which ensure the gluing of the (c) or (d) faces
with $\#A>2$ (see the proof of Lemma 5.4).

4) (see Proposition 6.7) If a connected component $A$ of a diagram $\g$ is
such that $A\cap U$ is an interval for $\s$, and if $\g'$ is a diagram
obtained from $\g$ by modifying $\s$ so that $A\cap U$ goes to another
cut of $\s$, then $\a (\g)=\a (\g')$.
\end{prop}

Now we are going to see that all the invariants produced by these 
constructions come from $Z_n$, provided that $\a $ is a map
from our set $\dl$ of diagrams to a free $\Z$-module and obeys the RSTU
relation. So the only interest in choosing such other maps $\a $ is
to obtain better rationality results than those of Theorem 
\ref{mult}.

\begin{prop} If a map $\a $ from $\dl$ to some $\R$-vector space
$\cal E$ obeys the RSTU relation, then the corresponding expression of
deg $\Phi$
is of the form $L(Z_n(K))$ where $L$ is a linear map from
$\A_n$ to $\cal E$ which only depends on $\a $.
\end{prop}

Let $\enk$ be the set of diagrams 
$\g\in\edn$ which are triply connected to $U$ and verify
$u_\g\geq k-1$.

We shall first define a map $\ps$ from $\enk$ to $\ce$. Then,
we shall prove that $\ps$ extends to a linear map $L$ defined
on $\A_n$, and that this $L$ is the solution to our problem.

Let us first define $\ps$ from $\enk$ to $\ce$ by:
$$\quad \ps(\g)=|\g|\mathop{\sum_{\g'\in\dl}}\limits_{\g'\sim \g}
\a (\g').$$

The existence of a linear extension of $\ps$ on $\A_n$ will be deduced
from the following facts:

1) Any RSTU relation implies the corresponding STU relation for $\ps$.

2) This map $\ps$ is well-defined on chord diagrams and obeys the
four-terms relation.

3) The STU relation obtained from 1) allows us to express any diagram 
in $\enk$ as a linear combination of chord diagrams.

Indeed, from 2) we deduce that $\ps$ defines a linear map from the space of
degree $n$ chord diagrams modulo the four-terms relation. But we know this
space is identified with $\A_n$, so $\ps$ defines a linear map $L$ on $\A_n$.
This map coincides with $\ps$ on $\enk$, thanks  to 1) and 3).

Now, to prove the identity
$$\sum_{\g\in \dl}\a (\g) I_K(\g)=\sum_{\g\in\edn}{I_K(\g)\over|\g|} L([\g])$$
we have to check that if a diagram $\g$ does not belong to $\enk$,
then $I_K(\g)=0$ or $L([\g])=0$.
Suppose $I_K(\g)\not=0$. Then $\g$ is triply connected to $U$, hence
$u_\g<k-1$. Thus, $\g$ can be expressed modulo STU 
as a linear combination of diagrams which have $k-1$ univalent 
vertices and are triply connected to $U$. But since $L$ vanishes on them,
it also vanishes on $\g$.
\medskip

Now let us prove the claims 1), 2) and 3).

{\it Proof of 1)}. Consider a RSTU relation; let $\g$ be any fixed diagram
in this relation, and let $A$ be the set of the two vertices involved.

The terms of the relation can be divided into three groups, which correspond to
the three terms of the STU relation.
Two of them consist of one single diagram (each of the two diagrams on the
left-hand side of the relation), and the third group is the sum on the
right-hand side. Let $T$ be the group which contains $\g$.
The diagram(s) in $T$ are exactly the diagram(s) $\g'$ such that 
$\g'$ is isomorphic to $\g$ by an isomorphism which fixes (with its
orientation) any visible edge which is not in $E_A$. They are
in one-to-one correspondence with the set of these isomorphisms.

Note that this RSTU relation is defined by the labelled diagram $\g/A$
(its diagrams are exactly the labelled diagrams which give $\g/A$ when
quotiented by an edge). We shall consider the sum of this relation
with respect to the relabellings of $\g/A$, in other words, with respect
to the $(\Upsilon,\lambda)$, where $\Upsilon$ is a diagram isomorphic to
$\g/A$, and $\lambda$ is an isomorphism from $\g/A$ to $\Upsilon$.

This sum of relations, when applied to $T$, is a sum over the set
of relabellings of $\g$, that is, the set of couples
$(\g',\lambda)$ where $\lambda$ is an isomorphism from $\g$ to $\g'$.
This set is in one-to-one correspondence with the product of the 
set of automorphisms
of $\g$ with the set of diagrams which are isomorphic to $\g$.
So, this sum applied to $T$ gives precisely $\ps(\g)$, by definition of
$\ps$. This proves 1).
\medskip

{\it Proofs of 2) and 3)}.
Since $k\leq 2n$, $\ps$ is well-defined on the chord diagrams, and also
on the diagrams with only one trivalent vertex. Then the four-terms relation
is a consequence of the STU relation deduced from 1).

The proof of 3) is easy by induction: one just has to check that when 
we use a STU relation to express a diagram in terms of two diagrams 
with one more 
univalent vertex, then these two diagrams are triply connected to 
$U$ provided the first one is.\cqfd
\medskip

This can be generalised to the case when $\a $ takes values in any free
$\Z$-module $M$, for $M$ is canonically imbedded in the real vector
space $M\otimes\R$.

\subsection{Rationality results in degree 2}

The present theorem gives the same rationality result as the one of 
\cite{le}, that is, a denominator of 48. But the actual denominator is 24,
according to the following formula deduced from \cite{bn'} (it uses
the value of the integral for the circle which was computed in \cite{gua}):
$$Z_2={I_K(\th)^2\over8}[\th]^2+(v_2+{1\over24})\left[\vep{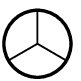}\right]$$
where $v_2 \in \Z$ is the degree 2 Vassiliev invariant which vanishes
on the unknot and takes the value 1 on the trefoils.

This result can be improved as did Polyak and Viro \cite{pv}.
We recall here their remarks, and precise the general context for which
this improvement in degree two is a particular case.

Keep $N=3$, but
consider the map $\a $ from $D_{2,3}$ into $\Z$ defined by:
$$\a \epfb{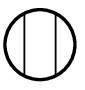}=0, \quad \a \epfb{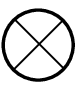}={1\over24}$$
independently of the labelling, and $$\a \epfb{p2}= 0\hbox{ or }-{1\over24}$$
according to the following rule: the labellings of this diagram are
divided into two classes according to the parity of the number
of (oriented) edges which point to the trivalent vertex. Then give the value
$-1\over24$ to a class and $0$ to the other.
It is easy to see that the gluing conditions (\ref{tot}) are realised for this
map $\a $ (see the weaker gluing condition given in the proof
of Lemma 5.4).

The expression of $Z_2$ suggests that it would be more relevant to study
the denominators of $Z(K)\over Z(O)$, where $O$ is the unknot.

\begin{lemma} The expression $$Z(K)\over Z(O)$$ is equal to the
expression defined the same way as $Z$, with a {\em long knot}:
instead of an imbedding of $S^1$ into $\rt$, it will be the imbedding
of $\R$ into $\rt$ which coincides with 
the inclusion of a given straight line $\lo$ of $\rt$ outside some compact set.
\end{lemma}

Let us sketch the proof of this lemma. Anyway, it will be a direct
consequence of the constructions in our next article.

Consider a knot which looks like a long knot when watched at some scale,
but is closed as a circle $O$ at some bigger scale. To prove that the integral
$Z$ of the whole is the product of the integrals of each ``part'', we
have to check that for a connected diagram $\g$ (see Remark 4.3 and 
Proposition 4.4), the configurations which 
take place at both scales do not contribute to the integral.

Let $A$ be the set of vertices which are mapped to the small scale. The
map $\Phi$ factorises through
$$C_{\g/A}(O)\lra (\sd)^{E^c_A},$$
and $\dim C_{\g/A}(O)\leq \dim(\sd)^{E^c_A}$. But in fact we are only
concerned with the configurations which map $A$ to a fixed point of $O$:
they form a codimension 1 part of $C_{\g/A}(O)$.
\cqfd

For this new shape of knot, we have to replace the cyclic order on the
univalent vertices of the labelled diagrams, by a linear order.
We are going to see the following lemma:

\begin{lemma} Rewrite the sufficient conditions 
(\ref{cond}), with the case of a linear order on the
univalent vertices, in a way which preserves this order.
Then the configuration spaces of a long knot 
can be glued, and $\Phi$ admits a degree.
\end{lemma}

Let us first define the compactification of the configuration spaces $\c$
we use here. We imbed it into the compact space 
$$\h\times[-\infty, \1]^U.$$

We have to check that the faces of this compactified space which correspond
to the $(f,f')\in \cb$ such that $f'(U)$ meets $\pm\infty$, 
are degenerate.
We suppose again (thanks to Remark 4.3) that the diagram $\g$ is connected.
Now, we deduce that $f_V$ maps all univalent vertices to the vertical 
straight line $\lo$. So, $f_V$ belongs to some space with dimension
$\leq 2\#E(V/S)-2$, which proves that the image by $\Phi$ of the 
corresponding stratum has codimension at least two.

Now, let us glue the anomaly faces. Note that the tangent map
d$K: [-\infty, \1]\lra S^2$ verifies
d$K(-\infty)=$ d$K(+\infty)$, so defines a map from
$S^1$ to $S^2$, which is the boundary of some map
$\ph: \dd\lra S^2$ in the same way as in 
Subsection \ref{gen}. With the same method, we can now prove that
$I_K(\th)+2I_\ph$ is an even integer.
\cqfd
\medskip

With this new problem, if we choose a map $\a $ which only depends
on the diagram with a cyclic order obtained by closing the linear
order, then this gives exactly the same rationality result as the one
we would obtain with this map from the diagrams with a cyclic
order, for the compact knot. So, we can obtain a better rationality result 
for $Z(K)\over Z(O)$ than for $Z(K)$, only when looking for maps $\a $ which do
depend on the cut of the cyclic order.

In particular for $n=2$ and $N=3$, take the map $\a $
which only depends on the orientations of the edges, such that
$$\a \epfb{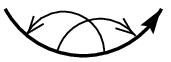}={1\over 6},\qquad \a \epfb{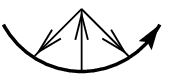}=-{1\over6}$$
and which vanishes on any other diagram. The corresponding expression
$\deg \Phi$ is equal to $v_2$.

We have verified that no better result can be obtained with another 
map $\a $ when $N=3$.

For more references on Vassiliev invariants, see Bar-Natan's bibliography:

http://www.ma.huji.ac.il/$\scriptstyle\sim$drorbn/VasBib/VasBib.html

\vfill
\begin{flushleft}
%adresse de l'auteur
{\footnotesize Sylvain POIRIER\\
INSTITUT FOURIER\\
Laboratoire de Math{\'e}matiques\\
UMR 5582 (UJF-CNRS)\\
BP 74\\
38402 St MARTIN D'H\`ERES Cedex (France)}
\end{flushleft}

\end{document}